\documentclass[12pt]{article}
\usepackage{amssymb, amsfonts, amsmath, amsthm, url, graphicx}

\def\3{\subset }
\def\4{\subseteq }

\def\cala{{\cal A}}
\def\0{\leqno}

\def\barr{\begin{array}}
\def\earr{\end{array}}
\def\dd{\displaystyle}

\def\Z{{\rlap{$\kern2pt{\rm Z}$}{\rm Z}\,}}

\title{\bf On the number of diamonds\\ in the subgroup lattice of a\\ finite abelian group}
\author{Dan Gregorian Fodor, Marius T\u arn\u auceanu}
\date{December 5, 2013}

\begin{document}

\maketitle

\begin{abstract}
The main goal of the current paper is to determine the total
number of diamonds in the subgroup lattice of a finite abelian
group. This counting problem is reduced to finite $p$-groups.
Explicit formulas are obtained in some particular cases.
\end{abstract}
\bigskip


\bigskip

\section{Introduction}

The relation between the structure of a group and the structure of
its lattice of subgroups constitutes an important domain of
research in group theory. The topic has enjoyed a rapid
development starting with the first half of the '20 century. Many
classes of groups determined by different properties of partially
ordered subsets of their subgroups (especially lattices of
subgroups) have been identified. We refer to Suzuki's book
\cite{11}, Schmidt's book \cite{7} or the more recent book
\cite{12} by the second author for more information about this
theory.
\bigskip

An important concept of subgroup lattice theory has been
introduced by Schmidt \cite{8} (see also \cite{9}): given a
lattice $L$, a group $G$ is said to be \textit{L-free} if the
subgroup lattice $L(G)$ has no sublattice isomorphic to $L$.
Interesting results about $L$-free groups have been obtained for
several particular lattices $L$, as the diamond lattice $M_5$ or
the pentagon lattice $N_5$. We recall here only that a group is
$M_5$-free if and only if it is locally cyclic and, in particular,
a finite group is $M_5$-free if and only if it is cyclic. Notice
also that the class of $L$-free groups can be extended to the
class of groups whose subgroup lattices contain a certain number
of sublattices isomorphic to $L$ (see e.g. \cite{15}).
\bigskip

Clearly, for a finite group $G$ the above concept leads to the
natural pro\-blem of counting the sublattices of $L(G)$ that are
isomorphic to a given lattice $L$. In the general case this
problem is very difficult. This is the reason for which in the
current paper we will treat only the particular case when $G$ is a
finite abelian group and $L=M_5$. In other words, we will
determine the number $dm(G)$ of all diamonds in the subgroup
lattice of a finite abelian group $G$. Recall that a \textit{diamond in $L(G)$}
(also called a \textit{diamond of $G$}) is a sublattice of $L(G)$
which is isomorphic to $M_5$. It is easy to see that such a sublattice
is completely determined by a triple $(A,B,C)\in L(G)^3$ satisfying
$A\wedge B=B\wedge C=C\wedge A$ and $A\vee B=B\vee C=C\vee A$.
\bigskip

The paper is organized as follows. In Section 2 we show that the
study can be reduced to $p$-groups and we develop a general method
to find $dm(G)$ for finite abelian $p$-groups. Section 3 deals
with the particular cases of ele\-men\-ta\-ry abelian $p$-groups
and of abelian $p$-groups of rank 2. The computation of $dm(G)$ is
also exemplified for a finite rank 4 abelian 2-group. In the final
section some conclusions and further research directions are
indicated.
\bigskip

Most of our notation is standard and will usually not be repeated
here. Basic definitions and results on lattices and groups can be
found in \cite{3,4} and \cite{6,10}, respectively.

\section{The key results}

Let $G$ be a finite abelian group and $L(G)$ be the subgroup
lattice of $G$. It is well-known that $G$ can be written as a
direct product of $p$-groups
$$G=\prod_{i=1}^k G_i,$$
where $|G_i|=p_i^{\alpha_i}$, for all $i=1,2,...,k$. Since the
subgroups of a direct product of groups having coprime orders are
also direct products (see Corollary of (4.19), \cite{10}, I), it
follows  that
$$L(G)\cong\prod_{i=1}^k L(G_i).$$The above lattice direct decomposition
is often used in order to reduce many combinatorial problems on
$L(G)$ to the subgroup lattices of finite $p$-groups (see e.g.
\cite{2,13,14}). This can be also applied to our problem.
\bigskip

First of all, let us assume that $k=2$. The following theorem
shows the way in which $dm(G)$ depends on $dm(G_1)$ and $dm(G_2)$.

\bigskip\noindent{\bf Theorem 2.1.} {\it Let $G$ be a finite abelian
group having a direct decomposition of type $G=G_1\times G_2$ with
$|G_i|=p_i^{\alpha_i}$ , $i=1,2$, and $p_1, p_2$ distinct primes.
Then the number of diamonds in the subgroup lattice of $G$
satisfies
$$dm(G)=dm(G_1)|L(G_2)|+dm(G_2)|L(G_1)|+6\, dm(G_1)\,dm(G_2).\0(1)$$}

\noindent{\bf Proof.} We easily infer that the diamonds
of $G$ are of one of the following three types:
\begin{description}
\item[\hspace{5mm}a)] $(A\times H,B\times H,C\times H)$, where $(A,B,C)$ is a diamond of $G_1$ and $H\in L(G_2)$;
\item[\hspace{5mm}b)] $(H\times A,H\times B,H\times C)$, where $(A,B,C)$ is a diamond of $G_2$ and $H\in L(G_1)$;
\item[\hspace{5mm}c)] $(A\times A',B\times B',C\times C')$, $(A\times A',B\times C',C\times B')$, $(A\times B',B\times A',C\times C')$, $(A\times B',B\times C',C\times A')$, $(A\times C',B\times A',C\times B')$, $(A\times C',B\times B',C\times A')$, where $(A,B,C)$ is a diamond of $G_1$ and $(A',B',C')$ is a diamond of $G_2$.
\end{description} Obviously, this leads to the desired equality (1).
\hfill\rule{1,5mm}{1,5mm}

\bigskip\noindent{\bf Example.} For the abelian group
$G=\mathbb{Z}_2^2\times \mathbb{Z}_3^2$ we have
$$dm(G)=dm(\mathbb{Z}_2^2)|L(\mathbb{Z}_3^2)|+dm(\mathbb{Z}_3^2)|L(\mathbb{Z}_2^2)|+6\,dm(\mathbb{Z}_2^2)\,dm(\mathbb{Z}_3^2)=$$
$$\hspace{-36mm}=1\cdot 6+4\cdot 5+6\cdot 1\cdot 4=50\,.$$

Even if we will not give an explicit formula, it is clear that the
above result can be extended for an arbitrary $k\geq 2$, and
consequently the computation of $dm(G)$ is reduced to the
computation of $dm(G_i)$, $i=1,2,...,k$. So, in the following we
will focus only on finite abelian $p$-groups.
\bigskip

In order to find $dm(G)$ for a finite abelian $p$-group $G$ we
need some auxiliary results. Recall first that a \textit{section}
of $G$ is a quotient of a subgroup of $G$. More precisely, every
section $S$ of $G$ is perfectly determined by a pair $(H,K)\in
L(G)^2$ such that $H\subseteq K$. We easily infer that $S$ is
isomorphic to a subgroup, or, dually, to a quotient of $G$. The
following lemma gives a way of counting the number of sections of
a fixed type in $G$ (notice that it can be also used together with
the Goursat's lemma to count the subgroups of $G$).

\bigskip\noindent{\bf Lemma 2.2.} {\it Let $G$ be a finite abelian
$p$-group and $S$ be a section of $G$. Then the number of sections
of $G$ that are isomorphic to $S$ is $$n_S(G)=\sum_{T\leq G,\,
T\cong S} |\,L(G/T)|\,.\0(2)$$Moreover, these sections can be
divided into $n$ classes, one for each subgroup $T_n \cong S$  of
$G$, and there is a bijection between the sections $K/H$ of a
class $m$ and the subgroups $R$ of $G/T_m$ such that if $K/H$
corresponds to $R$, then $H\cong R$ and $G/K\cong (G/T_m)/R$.}

\bigskip\noindent{\bf Proof.} The proof is based on a famous result
due to Baer (see Theorem 8.1.4 of \cite {7} or Theorem 4.2 of
\cite{11}) which states that every finite abelian group is
self-dual. In this way, if $\cala_G$ is the set of subgroup-quotient
group pairs $(A,G/A)$ of $G$, then one can define a bijection
$\delta_G :\cala_G\rightarrow\cala_G$, $\delta_G((A,G/A))=(A',G/A')$,
where $A\cong G/A'$ and $A'\cong G/A$.

We have to count all pairs $(H,K)\in L(G)^2$ such that $H\subseteq K$ and
$K/H\cong S$, or equivalently all pairs $(H,K/H)\in\cala_K$ such that $K/H\cong S$.
Since these pairs can be split into disjoint classes, one for every subgroup $K$ of
$G$, and each preselected bijection $\delta_K$ is defined on the class $\cala_K$ and preserves
it, we infer that counting the pairs $(H,K/H)\in\cala_K$ satisfying $K/H\cong S$ is
the same as counting the pairs $(H',K/H')\in\cala_K$ satisfying $H'\cong S$.
These can be obtained by counting for every $H'\cong S$ the number of $K$'s satisfying  $H'\subseteq K$.
For a fixed $H'$, the number of such $K$'s is $|L(G/H')|$. Thus
$$n_S(G)=\sum_{H'\leq G,\,H'\cong S} |\,L(G/H')|\,,$$which completes the proof.
\hfill\rule{1,5mm}{1,5mm}
\bigskip

Next we will focus on a particular class of diamonds in $L(G)$. A
diamond $(A,B,C)$ of $G$ will be called a \textit{primary diamond}
if $A\wedge B=B\wedge C=C\wedge A=1$ and $A\vee B=B\vee C=C\vee
A=G$. The structure of finite abelian $p$-groups having such a
diamond is very restrictive.

\bigskip\noindent{\bf Lemma 2.3.} {\it Let $G$ be a finite abelian
$p$-group and $(A,B,C)\in L(G)^3$ be a primary diamond. Then
$A\cong B\cong C\cong S$ and $G\cong S\times S$.}

\bigskip\noindent{\bf Proof.} It follows immediately that $G=A\times B=B\times C=C\times
A$. Then $A\cong G/C$ and $B\cong G/C$, and consequently $A\cong
B$. Similarly, one obtains $B\cong C$, completing the proof.
\hfill\rule{1,5mm}{1,5mm}
\bigskip

An explicit formula for the number of primary diamonds of a finite
abelian $p$-group $G$ of type $S\times S$ is indicated in the
following lemma.

\bigskip\noindent{\bf Lemma 2.4.} {\it Let $G\cong S\times S$ be a finite abelian
$p$-group. Then the number of primary diamonds in $L(G)$ is
$$\frac{|\,Aut(S\times S)|}{6\,|\,Aut(S)|}\,.\0(3)$$}

\bigskip\noindent{\bf Proof.} Write $G\cong S_1\times S_2$ with $S_1\cong S_2\cong S$.
We count first the pairs $(A,B)\in L(G)^2$ satisfying
$A\wedge B=1$ and $A\vee B=G$. It is clear that every automorphism
$f$ of $G$ determines such a pair $(A,B)$ given by $A=f(S_1)$ and
$B=f(S_2)$. Conversely, a pair $(A,B)$ corresponds to at least one
automorphism of $G$. However, many automorphisms can map to the
same such pair. More exactly, if two automorphisms $f_1$ and $f_2$
map to $(A,B)$, then the transformation between them is of type
$(g_1,g_2)\in Aut(S)^2$. Therefore the number of pairs $(A,B)\in
L(G)^2$ such that $A\wedge B=1$ and $A\vee B=G$ is
$$\frac{|\,Aut(S\times S)|}{|\,Aut(S)|^2}\,.$$

We still need to form triples $(A,B,C)\in L(G)^3$ with $A\wedge
B=B\wedge C=C\wedge A=1$ and $A\vee B=B\vee C=C\vee A=G$. Assume
next that the pair $(A,B)$ is fixed. Then $G\cong A\times B$ and
$A\cong B\cong S$. We easily infer that the number of $C$'s
completing $(A,B)$ to a primary diamond of $G$ is equal to the
number of isomorphisms between $A$ and $B$, or, equivalently, to
the number of automorphisms of $S$. Thus the number of the above
triples is
$$|\,Aut(S)|\,\,\frac{|\,Aut(S\times
S)|}{|\,Aut(S)|^2}=\frac{|\,Aut(S\times
S)|}{|\,Aut(S)|}\,.$$However, these triples are ordered and
consequently $G$ has
$$\frac{|\,Aut(S\times S)|}{6\,|\,Aut(S)|}$$primary diamonds, as desired.
\hfill\rule{1,5mm}{1,5mm}
\bigskip

We are now able to establish our main principle of counting the
diamonds of finite abelian $p$-groups.

\bigskip\noindent{\bf Theorem 2.5.} {\it The number of diamonds in the subgroup lattice
of a finite abelian $p$-group $G$ is
$$dm(G)=\sum_{S} n_{S\times S}(G)\, \frac{|\,Aut(S\times S)|}{6\,|\,Aut(S)|}\,,\0(4)$$where $S$ runs over
all types of sections $S$ of $G$ and the numbers $n_{S\times
S}(G)$ are given by {\rm (2)}.}

\bigskip\noindent{\bf Proof.} Let $(A,B,C)$ be a diamond of $G$. Obviously, by putting
$H=A\wedge B=B\wedge C=C\wedge A$ and $K=A\vee B=B\vee C=C\vee A$,
one obtains a pair $(H,K)$ of subgroups of $G$ such that
$H\subseteq K$, that is a section $S'$ of $G$. Moreover, $(A,B,C)$
is primary in $S'$ and so $S'$ is of type $S\times S$ (where
$A\cong B\cong C\cong S$) by Lemma 2.3. We infer that $dm(G)$ can
be obtained by counting all primary diamonds in all sections of
type $S\times S$ of $G$. Notice that the number of these sections
can be computed by using Lemma 2.2, while the number of primary
diamonds in such a section follows by (3). Hence (4) holds.
\hfill\rule{1,5mm}{1,5mm}
\bigskip

Finally, we recall the well-known formula for the number of
automorphisms of a finite abelian $p$-group (see e.g. \cite{1,5}),
which will be used in all particular cases in Section 3.

\bigskip\noindent{\bf Theorem 2.6.} {\it Let
$G=\prod_{i=1}^n \mathbb{Z}_{p^{\alpha_i}}$ be a finite abelian
$p$-group, where $1\leq\alpha_1\leq\alpha_2\leq...\leq\alpha_n$.
Then $$|Aut(G)|=\prod_{i=1}^n (p^{a_i}-p^{i-1})\prod_{u=1}^n
p^{\alpha_u(n-a_u)}\prod_{v=1}^n
p^{(\alpha_v-1)(n-b_v+1)}\,,$$where $$a_r=max\{s\mid
\alpha_s=\alpha_r\}\, \mbox{ and }\, b_r=min\{s\mid
\alpha_s=\alpha_r\}\,,\, r=1,2,...,n\,.$$In particular, we have
$$|Aut(\mathbb{Z}_p^n)|=p^{\frac{n(n-1)}{2}}\prod_{i=1}^n (p^i-1)\,.$$}

\section{Counting diamonds for certain finite abelian $p$-groups}

In the following let $p$ be a prime, $n$ be a positive integer and
$\mathbb{Z}_p^n$ be an elementary abelian $p$-group of rank $n$
(that is, a direct product of $n$ copies of $\mathbb{Z}_p$). We
recall first the well-known formula that gives the number
$a_{n,p}(i)$ of subgroups of order $p^i$ in $\mathbb{Z}_p^n$,
namely
$$a_{n,p}(i)=\frac{(p^n-1)\cdots (p-1)}{(p^i-1)\cdots
(p-1)(p^{n-i}-1)\cdots (p-1)}\,,\, i=0,1,...,n\,.$$Write $n=2m+r$
with $r\in\{0,1\}$. Then the $S\times S$-sections of
$\mathbb{Z}_p^n$ are of type $\mathbb{Z}_p^{2i}$, $i=1,2,...,m$.
For every $i$, the subgroups of $\mathbb{Z}_p^n$ isomorphic to
$\mathbb{Z}_p^{2i}$ are in fact all subgroups of order $p^{2i}$ of
$\mathbb{Z}_p^n$, and these have the same quotient
$\mathbb{Z}_p^{n-2i}$. Thus
$$n_{\mathbb{Z}_p^{2i}}(\mathbb{Z}_p^n)=a_{n,p}(2i)\,|L(\mathbb{Z}_p^{n-2i})|=a_{n,p}(2i)\sum_{j=0}^{n-2i}
a_{n-2i,p}(j)$$by Lemma 2.2. On the other hand, Theorem 2.6 shows
that
$$\frac{|\,Aut(\mathbb{Z}_p^{2i})|}{6\,|\,Aut(\mathbb{Z}_p^2i)|}=\frac{1}{6}\,p^{\frac{i(3i-1)}{2}}\prod_{k=i+1}^{2i}
(p^k-1)\,.$$Then (4) leads to the following result.

\bigskip\noindent{\bf Theorem 3.1.} {\it The number of diamonds in the subgroup lattice
of the finite elementary abelian $p$-group $\mathbb{Z}_p^n$ is
given by
$$dm(\mathbb{Z}_p^n)=\frac{1}{6}\,\sum_{i=1}^{[\frac{n}{2}]} p^{\frac{i(3i-1)}{2}} a_{n,p}(2i)\left(\sum_{j=0}^{n-2i}
a_{n-2i,p}(j)\right) \prod_{k=i+1}^{2i} (p^k-1)\,.\0(5)$$}
\bigskip

We exemplify the equality (5) for $n=4$.

\bigskip\noindent{\bf Example.} We have
$$dm(\mathbb{Z}_p^4)=\frac{p(p^2+p+1)(p^4-1)(p^5-p^4+p+3)}{6}$$and in particular
$$dm(\mathbb{Z}_2^4)=735\,.$$
\smallskip

Next we will focus on computing the number of diamonds for finite
rank 2 abelian $p$-groups, that is for groups of type
$\mathbb{Z}_{p^{\alpha_1}}\times\mathbb{Z}_{p^{\alpha_2}}$,
$1\leq\alpha_1\leq\alpha_2$. The $S\times S$-sections are in this
case of the form
$\mathbb{Z}_{p^i}^2=\mathbb{Z}_{p^i}\times\mathbb{Z}_{p^i}$,
$i=1,2,...,\alpha_1$. For every $i$, we can easily check that
$\mathbb{Z}_{p^{\alpha_1}}\times\mathbb{Z}_{p^{\alpha_2}}$ has a
unique subgroup isomorphic to $\mathbb{Z}_{p^i}^2$ and its
quotient is of type
$\mathbb{Z}_{p^{\alpha_1-i}}\times\mathbb{Z}_{p^{\alpha_2-i}}$. By
Theorem 3.3 of \cite{14}, we infer that
$$n_{\mathbb{Z}_{p^i}^2}(\mathbb{Z}_{p^{\alpha_1}}\times\mathbb{Z}_{p^{\alpha_2}})=|L(\mathbb{Z}_{p^{\alpha_1-i}}\times\mathbb{Z}_{p^{\alpha_2-i}})|=\frac{1}{(p{-}1)^2}\,f_p(\alpha_1-i,\alpha_2-i)\,,$$where
$$f_p(x_1,x_2)=(x_2{-}x_1{+}1)p^{x_1{+}2}{-}(x_2{-}x_1{-}1)p^{x_1{+}1}{-}(x_1{+}x_2{+}3)p{+}(x_1{+}x_2+1)\,,$$for
all $0\leq x_1\leq \alpha_1$ and $0\leq x_2\leq \alpha_2$.
Moreover, we have
$$\frac{|\,Aut(\mathbb{Z}_{p^i}^2)|}{6\,|\,Aut(\mathbb{Z}_{p^i})|}=\frac{1}{6}\,p^{3i-2}(p^2-1)$$by Theorem 2.6.
Hence the following result holds.

\bigskip\noindent{\bf Theorem 3.2.} {\it The number of diamonds in the subgroup lattice
of the finite rank {\rm 2} abelian $p$-group
$\mathbb{Z}_{p^{\alpha_1}}\times\mathbb{Z}_{p^{\alpha_2}}$,
$1\leq\alpha_1\leq\alpha_2$, is given by
$$dm(\mathbb{Z}_{p^{\alpha_1}}\times\mathbb{Z}_{p^{\alpha_2}})=\frac{p+1}{6(p-1)}\sum_{i=1}^{\alpha_1} p^{3i-2} f_p(\alpha_1-i,\alpha_2-i)\,,\0(6)$$
where the quantities $f_p(x_1,x_2)$ are indicated above.}
\bigskip

Obviously, by using a computer algebra program a precise
expression of
$dm(\mathbb{Z}_{p^{\alpha_1}}\times\mathbb{Z}_{p^{\alpha_2}})$ can
be obtained from (6). This counting is more facile in some
particular cases.

\bigskip\noindent{\bf Corollary 3.3.} {\it For every $n\in
\mathbb{N}^*$, we have:
\begin{itemize}
\item[{\rm a)}] $dm(\mathbb{Z}_p\times\mathbb{Z}_{p^n})=n\dd\binom{p+1}{3}${\rm ;}
\item[{\rm b)}] $dm(\mathbb{Z}_{2^n}\times\mathbb{Z}_{2^n})=\dd\frac{1}{49}\left(3\cdot 2^{3n+2}-49\cdot 2^n+14n+37\right)$\,.
\end{itemize}}
\bigskip

We also exemplify the equality (6) for $\alpha_1=2$ and
$\alpha_2=3$.

\bigskip\noindent{\bf Example.} We have
$$dm(\mathbb{Z}_{p^2}\times\mathbb{Z}_{p^3})=\frac{p(p-1)(p+1)^2(p^2-p+2)}{3}$$and in particular
$$dm(\mathbb{Z}_4\times\mathbb{Z}_8)=24\,.$$
\smallskip

In the end of this section we will apply Theorem 2.5 to count the
diamonds of a rank 4 abelian 2-group, namely
$G=\mathbb{Z}_2\times\mathbb{Z}_4^3$. These can be divided into
four (disjoint) classes:
\begin{itemize}
\item[{\rm a)}] Those that are primary in $\mathbb{Z}_2\times\mathbb{Z}_2$-sections, 1 per section.
We can easily see that $G$ has 35 subgroups isomorphic with $\mathbb{Z}_2\times\mathbb{Z}_2$\,: 28 have a quotient of type
$\mathbb{Z}_2\times\mathbb{Z}_4^2$  and 7 have a quotient of type  $\mathbb{Z}_2^3\times\mathbb{Z}_4$.
Then Lemma 2.2 implies that the number of sections of type $\mathbb{Z}_2\times\mathbb{Z}_2$ in $G$ is
$$n_{\mathbb{Z}_2\times\mathbb{Z}_2}(G)=28\cdot|L(\mathbb{Z}_2\times\mathbb{Z}_4^2)|+7\cdot|L(\mathbb{Z}_2^3\times\mathbb{Z}_4)|=28\cdot54+7\cdot118=2338\,,$$while
the total number of primary diamonds in
$\mathbb{Z}_2\times\mathbb{Z}_2$-sections is $$2338\cdot
1=2338\,.$$
\item[{\rm b)}] Those that are primary in $\mathbb{Z}_4\times\mathbb{Z}_4$-sections,
$\frac{|Aut(\mathbb{Z}_4\times\mathbb{Z}_4)|}{6\,|Aut(\mathbb{Z}_4)|}=8$
per section. $G$ has 112 subgroups isomorphic with
$\mathbb{Z}_4\times\mathbb{Z}_4$, all of them having quotients of
type $\mathbb{Z}_2\times\mathbb{Z}_4$. It follows that that the
number of sections of type $\mathbb{Z}_4\times\mathbb{Z}_4$ in $G$
is
$$n_{\mathbb{Z}_4\times\mathbb{Z}_4}(G)=112\cdot|L(\mathbb{Z}_2\times\mathbb{Z}_4)|=112\cdot8=896\,,$$while
the total number of primary diamonds in
$\mathbb{Z}_4\times\mathbb{Z}_4$-sections is
$$896\cdot8=7168\,.$$
\item[{\rm c)}] Those that are primary in $\mathbb{Z}_2^4$-sections, $\frac{|Aut(\mathbb{Z}_2^4)|}{6\,|Aut(\mathbb{Z}_2^2)|}=560$ per section.
$G$ has 1 subgroup isomorphic with $\mathbb{Z}_2^4$, whose
quotient is isomorphic to $\mathbb{Z}_2^3$. It follows that the
number of sections of type $\mathbb{Z}_2^4$ in $G$ is
$$n_{\mathbb{Z}_2^4}(G)=1\cdot|L(\mathbb{Z}_2^3))|=1\cdot16=16\,,$$while the total number
of primary diamonds in $\mathbb{Z}_2^4$-sections is
$$16\cdot560=8960\,.$$
\item[{\rm d)}] Those that are primary in $\mathbb{Z}_2^2\times\mathbb{Z}_4^2$-sections, $\frac{|Aut(\mathbb{Z}_2^2\times\mathbb{Z}_4^2)|}{6\,|Aut(\mathbb{Z}_2\times\mathbb{Z}_4)|}=3072$ per section.
$G$ has 7 subgroups isomorphic with
$\mathbb{Z}_2^2\times\mathbb{Z}_4^2$, all of them having quotients
of type $\mathbb{Z}_2$. So, the number of sections of type
$\mathbb{Z}_2^2\times\mathbb{Z}_4^2$ in $G$ is
$$n_{\mathbb{Z}_2^2\times\mathbb{Z}_4^2}(G)=7\cdot|L(\mathbb{Z}_2)|=7\cdot2=14\,,$$while the total number of
primary diamonds in $\mathbb{Z}_2^2\times\mathbb{Z}_4^2$-sections
is $$14\cdot3072=43008\,.$$
\end{itemize}

Hence
$$dm(\mathbb{Z}_2\times\mathbb{Z}_4^3)=2338+7168+8960+43008=61474\,.$$

\section{Conclusions and further research}

All our previous results show that the problem of counting the
number of sublattices in the subgroup lattice of a group $G$ that
are isomorphic to a given lattice $L$ is an interesting
computational aspect of subgroup lattice theory. Clearly, the
study started in this paper can be extended for other lattices $L$
and groups $G$. This will surely be the subject of some further
research.
\newpage

Finally, we indicate three open problems concerning the above topic.
\bigskip

\noindent{\bf Problem 5.1.} Improve the results of Section 3, by
obtaining an explicit formula for $dm(G)$ when $G$ is a finite
abelian $p$-group of an arbitrary rank.
\bigskip

\noindent{\bf Problem 5.2.} In the class of finite groups $G$ of a
fixed order, find the minimum/maximum of $dm(G)$. Is it true that
the function $dm$ is strictly de\-crea\-sing on the set of types
abelian $p$-groups of order $p^n$, totally ordered by the
lexicographic order?
\bigskip

\noindent{\bf Problem 5.3.} Determine the number of sublattices in
other remarkable posets of subgroups (e.g. normal subgroup
lattices) of a finite group that are isomorphic to a given
lattice.
\bigskip

\bigskip\noindent{\bf Acknowledgements.} The authors are grateful to the reviewers for
their remarks which improve the previous version of the paper.

\vspace*{5ex} \small

\begin{minipage}[t]{6cm}
Dan Gregorian Fodor \\
Faculty of Mathematics \\
``Al.I. Cuza'' University \\
Ia\c si, Romania \\
e-mail: {\tt dan.ms.chaos@gmail.com}
\end{minipage}
\hfill
\begin{minipage}[t]{5cm}
Marius T\u arn\u auceanu \\
Faculty of  Mathematics \\
``Al.I. Cuza'' University \\
Ia\c si, Romania \\
e-mail: {\tt tarnauc@uaic.ro}
\end{minipage}

\end{document}